\documentclass[11pt]{amsart}
\usepackage{amsmath,amsfonts,amssymb,latexsym,amsthm}
\usepackage{epsfig}
\usepackage{hyperref} 
\usepackage[dvipsnames]{xcolor}

\newtheorem{thm}{Theorem}

\newtheorem{conv}[thm]{Convention}

\title{Positive dependence for colored percolation}

\thanks{\today}

\author[Nikita Gladkov \and Igor Pak]{Nikita Gladkov$^\ast$ \. and \. Igor Pak$^\ast$}

\thanks{\thinspace ${\hspace{-.99ex}}^\ast$Department of Mathematics,
 University of California, Los Angeles. \\ \ {\hskip.5cm}  Email: \ts \texttt{{\{gladkovna,pak\}@math.ucla.edu}}}

\def\lra{\leftrightarrow}

\def\emp{\varnothing}

\def\zz{\mathbb Z}

\def\sm{\smallsetminus}

\def\al{\alpha}

\def\cB{\mathcal B}

\def\cA{\mathcal A}

\def\cU{\mathcal U}
\def\cV{\mathcal V}

\def\cW{\mathcal W}

\def\bP{\mathbb{P}}

\def\<{\langle}
\def\>{\rangle}

\def\0{{\mathbf 0}}

\def\.{\hskip.06cm}
\def\ts{\hskip.03cm}

\def\mX{{\mathcal{X}}}

\begin{document}





\begin{abstract}
For uniform random $4$-colorings of graph edges with colors $\{a,b,c,d\}$,
every two colors form a $\frac{1}{2}$-percolation, and every two overlapping
pairs of colors form independent $\frac{1}{2}$-percolations.  We show joint
positive dependence for pairs of colors $ab$, $ac$ and~$ac$, and
joint negative dependence for pairs of colors $ab$, $ac$ and~$bc$.
The proof is based on a generalization of the Harris--Kleitman
inequalities.  We apply the results to crossing probabilities for
the colored bond and site percolation, and to colored critical
percolation that we also define.
\end{abstract}

\maketitle


\subsection*{Introduction}
The study of \emph{percolation} goes back to the 1957 paper by
Broadbent and Hammersley~\cite{BH57} and has been incredibly
popular in the last few decades across the sciences.  It remains
one of the most applied statistical models, reaching far corners
of statistical physics and probability, and fields as disparate
a materials science, network theory and seismology,
see e.g.\ \cite{Gri,SCA,Sah}.

Despite remarkable recent advances, many problems remain open
and continued to be actively pursued, see e.g.\ \cite{AGKSZ,Dum,Gri23,Mor}.
Note that specific models of percolation wary greatly depending on the
scientific context and applications.  Here we consider the
\emph{colored} bond $($site$)$ percolation, where each graph edge (vertex) takes
random color, see e.g.~\cite{KM,SCA,Zal}.

As one studies random events, one is naturally concerned about
their correlations.  The \emph{Fortuin--Kasteleyn--Ginibre {\rm (FKG)} inequality}
\cite{FKG71} is a basic tool to
establish positive dependence for percolation and related models,
see e.g.\ \cite{DS,New,Wer}.  This inequality shows that every two
increasing (or two decreasing) random events are positively correlated (see below).

The FKG inequality is itself an advance generalization
of the \emph{Harris--Kleitman inequality} \cite{Har,Kle66} discovered independently
in probability and graph theory.  Outside of its fundamental applications to
statistical physics and probability, it has numerous applications in graph
theory~\cite{Cha,JLR}, order theory \cite{Fish,She-XYZ} and algebraic
combinatorics~\cite{CP23}.

We are interested in generalizations of the Harris--Kleitman inequality
to multiple functions, which has also been intensely studied but remains largely
mysterious \cite{Gla,LS,Sahi}.  More precisely, we establish positive
dependence for three pairwise independent percolations and generalize it
further to \ts $k$ \ts percolations such that every \ts $(k-1)$ \ts
of them are mutually independent.

\smallskip

\subsection*{Positive correlation in percolation}
We first illustrate the power of the Harris--Kleitman inequality. Let $G = (V, E)$
be a simple graph, which can be finite or infinite. Consider a
\emph{$p$-percolation} defined by
independently at random deleting edges of~$G$ with probability $(1-p)$.
We write $\bP_p(x\lra y)$ for the probability that vertices $x,y\in V$
are connected.

In its basic application, the Harris--Kleitman inequality proves a positive correlation
of connectivity of two pairs of vertices:
\begin{equation}\label{eq:lra}
\bP_p(x \lra  y, \ts u \lra v) \, \ge \, \bP_p(x \lra y) \. \bP_p(u \lra v),
\end{equation}
for all $x,y,u,v\in V$.  Equivalently, this says that the probability
that two vertices are connected increases if \emph{some other} two
vertices are connected, even if these two vertices are quite far in
the graph: \. $\bP_p(x \lra  y \. | \. u \lra v) \. \ge \. \bP_p(x \lra y)$.
This easily implies that the \emph{critical probability} \ts
$p_c= \sup\big\{p \. : \. \bP_p(x\lra \infty)=0\big\}$ \ts is independent
on the vertex~$x$ in every connected graph, see e.g.\ \cite{BR,Gri}.  For
the case when $G = \zz^2$ is a square lattice, Harris proved that \ts
$p_c\ge \frac12$ \ts in the original paper~\cite{Har}.  Famously,
Kesten~\cite{Kes} established the equality \ts $p_c= \frac12$ \ts
twenty years later.

More generally, a subset $\cA \subseteq 2^E$ is called \emph{closed upward},
if $A+e\in \cA$ for every $A\in \cA$ and $e\in E\sm A$.  Similarly, $\cA$ is
\emph{closed downward}, if $A-e\in \cA$ for every $A\in \cA$ and $e\in A$.
We think of $\cA$ as \emph{graph
property}, and write $\bP_p(\cA)$ for the probability that the property holds for
a $p$-percolation.  In this notation, the \emph{Harris--Kleitman inequality}
 states that \.
\begin{equation}\label{eq:HK}
\bP_p(\cA \cap \cB) \ge \bP_p(\cA) \. \bP_p(\cB),
\end{equation}
for every
two closed upward \ts $\cA,\cB$.  For \ts $\cA=\{H \ts : \ts x\lra y\}$ \ts and
\ts $\cB=\{H \ts : \ts u\lra v\}$ \ts we obtain~\eqref{eq:lra}.
Note that \eqref{eq:HK} holds also for every two closed downward \ts $\cA,\cB$.

\smallskip

\subsection*{Positive dependence in colored percolation}
Let \ts $f: E\to \{a,b,c,d\}$ \ts be a uniform random coloring of the edges of~$G$,
where each edge is colored uniformly and independently. This gives a random partition \ts
$E = E_a \sqcup E_b \sqcup E_c \sqcup E_d$, where \ts $E_s=f^{-1}(s)$ \ts for a color \ts
$s\in \{a,b,c,d\}$.

For every two distinct colors \ts $s,t\in \{a,b,c,d\}$, denote \. $E_{st}:=E_s\cup E_t$\ts.
One can think of $E_{st}$ as either a $\frac12$-percolation or a random uniform spanning
subgraph of~$G$.  Note that $E_{ab}$, $E_{ac}$ and $E_{bc}$ are pairwise independent,
but not mutually independent.  Our main result establishes their negative dependence:

\smallskip

\begin{thm}\label{t:events}
Let \ts $\cU, \cV, \cW$ \ts be closed upward graph properties.  Denote by
\ts $\cU_{ab}$, \ts $\cV_{ac}$  \ts and \ts $\cW_{bc}$ \ts  the corresponding
properties of \ts $E_{ab}$, $E_{ac}$ and $E_{bc}$, respectively.  Then the events \ts $\cU_{ab}$,
\ts $\cV_{ac}$ \ts and \ts $\cW_{bc}$ \ts are pairwise independent, but
\begin{equation}\label{eq:main}
\bP(\cU_{ab} \cap \cV_{ac} \cap \cW_{bc}) \, \le \, \bP(\cU_{ab}) \. \bP(\cV_{ac}) \. \bP(\cW_{bc}),
\end{equation}
where the probability is over uniform random colorings \ts $f: E\to \{a,b,c,d\}$.  Similarly,
\begin{equation}\label{eq:main2}
\bP(\cU_{ab} \cap \cV_{ac} \cap \cW_{ad}) \, \ge \, \bP(\cU_{ab}) \. \bP(\cV_{ac}) \. \bP(\cW_{ad}),
\end{equation}
where \ts $\cW_{ad}$ \ts is the property of \ts $E_{ad}$\ts.
\end{thm}

\smallskip

Since all \ts $E_{st}$ \ts are $\frac12$-percolations, we can rewrite the RHS of both~\eqref{eq:main} and~\eqref{eq:main2}
as a more symmetric product:
\begin{equation}\label{eq:main-RHS}
\bP_{\frac12}(\cU) \. \bP_{\frac12}(\cV) \. \bP_{\frac12}(\cW).
\end{equation}

For example, let
$E=\{e\}$, so that $G$ is a graph with one edge, and let \ts $\cU=\cV=\cW$ \ts be properties
of containing~$e$.  The LHS of~\eqref{eq:main} is zero since we always have \ts
$E_{ab}\cap E_{ac} \cap E_{bc}=\emp$.  Similarly, the LHS of~\eqref{eq:main2} is \ts $\frac14$
\ts since $E_{ab}\cap E_{ac} \cap E_{ad}=E_a$.  On the other hand, the product~\eqref{eq:main-RHS}
 is \ts $\frac{1}{8}$ \ts
since \ts $\bP_{\frac12}(\cU) = \bP_{\frac12}(\cV) = \bP_{\frac12}(\cW) = \frac12$\ts.

\smallskip

\subsection*{Proof of Theorem~\ref{t:events}.}
Since \ts $E_{ab}$ \ts and \ts $E_{ac}$ \ts are independent $\frac12$-percolations,
this implies that events $\cU_{ab}$ and $\cV_{ac}$ are also independent.  This
proves the pairwise independence part.

We prove \eqref{eq:main} by induction on the number of edges in $E$. For $E=\emp$, the inequality is trivial.
Fix an edge $e \in E$. Consider the probability space of colorings of $E-e$. For an event \ts
$\mX_{ab} \subseteq 2^E$,
denote by \ts $\mX^+_{ab}$ \ts the subset of \ts
$\mX_{ab}$ \ts such that \ts $f(e)\in \{a,b\}$.  Similarly,
denote by \ts $\mX^-_{ab}$ \ts the subset of \ts
$\mX_{ab}$ \ts such that \ts $f(e)\in \{c,d\}$.

By the symmetry, we have: \.
$$\bP\big(\mX_{ab}\.:\. f(e)=a\big) \. = \. \bP(\mX_{ab}\.:\. f(e)=b) \. = \.  2\ts\bP_\frac12(\mX^+),
$$
$$\bP\big(\mX_{ab}\.:\. f(e)=c\big) \. = \. \bP(\mX_{ab}\.:\. f(e)=d) \. = \.  2\ts\bP_\frac12(\mX^-).
$$
Clearly, \. $\bP_\frac12(\mX)= \bP_\frac12(\mX^-) + \bP_\frac12(\mX^+)$.
When \ts $\mX$ \ts is closed upward, we also have \ts $\bP_\frac12(\mX^-) \le \bP_\frac12(\mX^+)$.
We use this notation for \ts $\mX\in \{\cU,\cV,\cW\}$ \ts and all pairs of colors.

Considering all possible colors of $e$ and using the induction hypothesis, we have:
\begin{align*}
\bP(\cU_{ab} \cap \cV_{ac} \cap \cW_{bc}) \,
& = \,  \ts \bP\big(\cU_{ab}^+ \cap \cV_{ac}^+ \cap \cW_{bc}^-\big) \. + \.  \ts  \bP\big(\cU_{ab}^+ \cap \cV_{ac}^- \cap \cW_{bc}^+\big) \\
&\hspace{6em}+  \ts  \bP\big(\cU_{ab}^- \cap \cV_{ac}^+ \cap \cW_{bc}^+\big) \. + \.   \ts \bP\big(\cU_{ab}^- \cap \cV_{ac}^- \cap \cW_{bc}^-\big) \\
&\le \,  2 \Big(\ts \bP\big(\cU_{ab}^+ \big) \. \bP\big( \cV_{ac}^+ \big) \. \bP\big( \cW_{bc}^-\big) \. +   \ts \bP\big(\cU_{ab}^+\big) \. \bP\big( \cV_{ac}^-\big) \. \bP\big( \cW_{bc}^+\big) \\
&\hspace{6em}+   \ts \bP\big(\cU_{ab}^- \big) \.\bP\big( \cV_{ac}^+\big) \.\bP\big( \cW_{bc}^+\big) \. + \.   \ts\bP\big(\cU_{ab}^-\big) \. \bP\big( \cV_{ac}^- \big) \bP\big( \cW_{bc}^-\big)\Big).
\end{align*}
Simplifying the notation as above, the RHS is equal to:
\begin{align*}
& 2 \Big(\ts \bP_\frac12(\cU^+) \. \bP_\frac12( \cV^+ ) \. \bP_\frac12( \cW^-) \. +  \ts \bP_\frac12(\cU^+) \. \bP_\frac12( \cV^-) \. \bP_\frac12( \cW^+) \\
&\hspace{8em}+  \ts \bP_\frac12(\cU^-) \.\bP_\frac12( \cV^+) \.\bP_\frac12( \cW^+) \. + \.  \ts\bP_\frac12(\cU^-) \. \bP_\frac12( \cV^- ) \bP_\frac12( \cW^-) \Big)\\
&\hspace{1em} = \, \big(\bP_\frac12(\cU^+)+\bP_\frac12(\cU^-)\big)\big(\bP_\frac12(\cV^+)+\bP_\frac12(\cV^-)\big)\big(\bP_\frac12(\cW^+) + \bP_\frac12(\cW^-)\big)\\
&\hspace{8em}-\big(\bP_\frac12(\cU^+)-\bP_\frac12(\cU^-)\big)\big(\bP_\frac12(\cV^+) - \bP_\frac12(\cV^-)\big)\big(\bP_\frac12(\cW^+)  - \bP_\frac12(\cW^-)\big) \\
&\hspace{1em}\le \, \bP_\frac12(\cU) \. \bP_\frac12(\cV) \. \bP_\frac12(\cW),
\end{align*}
as desired.  The proof of~\eqref{eq:main2} goes along the same lines.\qed

\smallskip

\subsection*{Variations and generalizations}
First, note that we never use the graph structure, and the theorem can be viewed as a
result about abstract set systems, cf.~\cite{AS,Kle66}.  In particular, it applies
to the site percolation (see below).
Note also that the theorem can be extended to the $p$-percolation for all \ts $0 \le p\le 1$,
but the resulting coupling of percolations then require seven colors
and have somewhat inelegant probabilities \cite{Gla-thesis}.

Curiously, for \emph{closed downward} \ts properties, the inequalities in the theorem hold reverse:
\begin{equation}\label{eq:main3}
\bP(\cU_{ab} \cap \cV_{ac} \cap \cW_{bc}) \, \ge \, \bP(\cU_{ab}) \. \bP(\cV_{ac}) \. \bP(\cW_{bc})
\end{equation}
and
\begin{equation}\label{eq:main4}
\bP(\cU_{ab} \cap \cV_{ac} \cap \cW_{ad}) \, \le \, \bP(\cU_{ab}) \. \bP(\cV_{ac}) \. \bP(\cW_{ad}).
\end{equation}
The proofs follow verbatim the proofs in the theorem.

\smallskip

Next, we generalize the theorem to larger number of events.
Start by taking $k$ independent $\frac12$-percolations  \ts $E_1, \dots, E_k$ \ts on the same graph.
Define a new \ts $\frac12$-percolation \. $E_{k+1} := \bigoplus E_i\mod 2$,
where the edge $e$ is present if and only if it is present in an odd number of $E_i$'s.
Observe that every $k$ of \ts $E_1,\ldots,E_{k+1}$ \ts are mutually independent.

Then, for every closed downward properties \. $\mX_1, \dots, \mX_{k+1}$ \. we have:
\begin{equation}\label{eq:main-multi}
\bP\big(\mX_1 \cap \dots \cap \mX_{k+1}\big) \, \ge \, \bP(\mX_1) \. \cdots \. \bP(\mX_{k+1}).
\end{equation}
Once again, the proof follows verbatim the proof of the theorem.

Note that for $k=1$, we have \ts $E_1=E_2$ \ts and \eqref{eq:main-multi}
is the Harris--Kleitman inequality~\eqref{eq:HK}.  For \ts $k=2$, let \ts
$$f(e) \ := \ \left\{
\aligned
& \, a \quad \text{if} \ \ e\in E_1 \cap E_2 \\
& \, b \quad \text{if} \ \ e\in E_1, \. e\notin E_2 \\
& \, c \quad \text{if} \ \ e\in E_2, \. e\notin E_1 \\
& \, d \quad \text{if} \ \ e\notin E_1, \. e\notin E_2
\endaligned \right.
$$
Taking \ts $\cU:=\mX_1$\ts, \ts $\cV:=\mX_2$ \ts and \ts $\cW:=\mX_3$\ts,
we have~\eqref{eq:main-multi} coincides with~\eqref{eq:main3}.

Finally, one can easily obtain a colored version with $m=2^k$ colors.
E.g., for \ts $k=3$, take a uniform random coloring \ts $f:E\to \{1,\ldots,8\}$.  Consider
four pairwise independent \ts $\frac12$-percolations \ts $E_{1234}$,
\ts $E_{1256}$, \ts $E_{1357}$ \ts and \ts $E_{1467}$ \ts with natural
labeling.  Note that every three of these are mutually independent.
Then, for closed downward properties \ts $\cU,\cV,\cW$ \ts
and \ts $\mX$,  the inequality~\eqref{eq:main-multi} gives:
\begin{equation*}
\bP(\cU_{1234} \cap \cV_{1256} \cap \cW_{1357} \cap \mX_{1467}) \,
\ge \, \bP(\cU_{1234}) \. \bP(\cV_{1256}) \. \bP(\cW_{1357}) \. \bP(\mX_{1467}).
\end{equation*}

\smallskip

\subsection*{Probability of the majority}
The simplest nontrivial example in the theorem is when \ts $\cU=\cV=\cW$ \ts
is the property of having \ts $>m$ \ts edges, where \ts $|E|=2m+1$.  The graph structure
is irrelevant in this case, and we have \ts $\bP_{\frac12}(\cU_{ab})=\frac12$\ts.
A direct calculation in this case gives:
\begin{equation*}
\bP(\cU_{ab} \cap \cV_{ac} \cap \cW_{bc}), \,  \bP(\cU_{ab} \cap \cV_{ac} \cap \cW_{ad}) \,
\to \, \bP_{\frac12}(\cU)^3 \, = \, \tfrac{1}{8}
\end{equation*}
as \ts $m\to \infty$.  This shows that both \eqref{eq:main} and \eqref{eq:main2}
are asymptotically tight in this case.

\smallskip

\subsection*{Crossing probabilities in a rectangle}
Let \ts $G=(V,E)$ \ts be a $n\times (n+1)$ rectangle as in Figure~1.
Consider a uniform random coloring \ts $f: E\to \{a,b,c,d\}$.
Note that \ts $E_{ab}$, \ts $E_{ac}$ \ts and \ts $E_{ad}$ \ts are pairwise
independent bond $\frac{1}{2}$-percolations with free boundary conditions~(BC).
Let \ts $\cU=\{12\lra 34\}$ \ts be the connectivity
property of the opposite sides of~$G$, and recall that \ts
$\bP_{\frac12}(\cU_{ab})=\frac12\ts$, see e.g.\ \cite{BR}.  Then \eqref{eq:main2} gives:
\begin{equation}\label{eq:square-cross}
\bP\big(\cU_{ab} \cap \cU_{ac} \cap \cU_{ad}\big) \, \ge \, \bP_{\frac12}(\cU)^3 \, = \, \tfrac18\ts,
\end{equation}
for all~$n\ge 1$.  On the other hand, by the pairwise independence we have:
$$
\bP\big(\cU_{ab} \cap \cU_{ac} \cap \cU_{ad}\big) \, \le \,\bP\big(\cU_{ab} \cap \cU_{ac}\big)
\, = \, \bP_{\frac12}(\cU)^2 \, = \, \tfrac14\ts.
$$
Note that as a function of~$p$ the crossing probability in a rhombus under
\ts $p$-percolation has a sharp threshold \cite{BR}, so the trivial lower bound
is unhelpful:
$$
\bP\big(\cU_{ab} \cap \cU_{ac} \cap \cU_{ad}\big) \. \ge \. \bP(\cU_a) \. = \. \bP_{\frac14}(\cU)
\. \to \. 0 \quad \text{as \ $n\to \infty$.}
$$
For \ts $n=30$, the sampling of \ts $N=4\cdot 10^7$ \ts trials
gives an approximation \. $\bP\big(\cU_{ab} \cap \cU_{ac} \cap \cU_{ad}\big) = 0.125098 \pm 0.000052$.
We conjecture that this probability is $\ts\frac18$ \ts in the limit \ts $n\to \infty$.

\begin{figure*}[hbt]
\begin{center}
	\includegraphics[height=2.1cm]{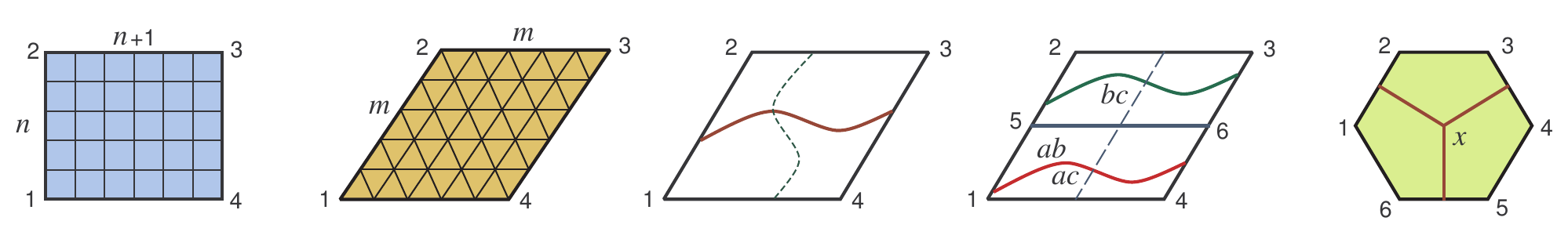}
\caption{Crossing probabilities in a rectangle, rhombus and a hexagon.}
\end{center}
\label{f:tri}
\end{figure*}

\smallskip

\subsection*{Crossing probabilities in a rhombus}
Let \ts $G=(V,E)$ \ts be a $m$-rhombus on the triangular
lattice, see Figure~1.  Consider a uniform random coloring \ts $f: V\to \{a,b,c,d\}$.
Note that \ts $V_{ab}$, \ts $V_{ac}$ \ts and \ts $V_{ad}$ \ts are pairwise
independent site $\frac{1}{2}$-percolations with free BC.
Let \ts $\cU=\{12\lra 34\}$ \ts and  \ts $\cU'=\{14\lra 23\}$ \ts be connectivity
properties of the opposite sides of~$G$.  Recall that \ts
$\bP_{\frac12}(\cU_{ab})+\bP_{\frac12}(\cU'_{cd})=1$ \ts by a topological argument,
so \ts $\bP_{\frac12}(\cU)=\bP_{\frac12}(\cU')=\frac12$ \ts by the symmetry.  Then
\eqref{eq:main} and \eqref{eq:main2} give:
\begin{equation}\label{eq:square-cross}
\aligned
\bP\big(\cU_{ab} \cap \cU_{ac} \cap \cU_{bc}\big) \, & \le \, \bP_{\frac12}(\cU)^3 \, = \, \tfrac18\ts,\\
\bP\big(\cU_{ab} \cap \cU_{ac} \cap \cU_{ad}\big) \, & \ge \, \bP_{\frac12}(\cU)^3 \, = \, \tfrac18\ts,
\endaligned
\end{equation}
for all~$m\ge 1$.
We conjecture that
$$\bP\big(\cU_{ab} \cap \cU_{ac} \cap \cU_{bc}\big) \. \to \. \tfrac18 \quad \text{and} \quad
 \bP\big(\cU_{ab} \cap \cU_{ac} \cap \cU_{ad}\big) \. \to \. \tfrac18 \quad \text{as  \ $m\to \infty$.}
$$
If this holds, we also have other similar limits, e.g.
$$
\bP\big(\cU_{ab} \cap \cU_{ac} \cap \cU'_{bc}\big) \,  = \, \bP_{\frac12}(\cU)^2 \. - \. \bP\big(\cU_{ab} \cap \cU_{ac} \cap \cU_{ad}\big) \, \to \, \tfrac18\..
$$
This is in contrast with limits such as \. $\bP\big(\cU_{ab} \cap \cU_{bc} \cap \cU_{cd}\big)$ \.
which can be computed by \emph{Watts' formula} \cite{Watts} (see also~\cite{Dub,SW}).
Finally, we note that \ts
$\bP\big(\cU_{ab} \cap \cU_{ac} \cap \cU_{bc}\big)$ \ts is bounded away from zero.
To see this, partition the  rhombus into four parallelograms (see Figure~1), so the desired
probability is bounded by the crossing probabilities:
$$\bP\big(\cU_{ab} \cap \cU_{ac} \cap \cU_{bc}\big) \. \ge \.
2 \. \bP_{\frac12}\big(15 \lra_{ab} 46, \ts 15 \lra_{ac} 46\big) \. \bP_{\frac12}\big(52 \lra_{bc} 63\big)
\, \ge  \, 2 \. \tfrac1{4^2} \. \tfrac1{2^2} \. = \. \tfrac1{32}\..
$$
This bound can be improved to \ts $\frac{9}{128}$ \ts by a careful use of the inclusion-exclusion.
In the limit \ts $m\to \infty$, this bound can be further improved since
these crossing probabilities can be computed by \emph{Cardy's formula} \cite{BR,Gri}.

\smallskip

\subsection*{Crossing probabilities in a hexagon}
Consider a regular hexagon \ts $G=(V,E)$ \ts on the triangular lattice with
side lengths~$\ell$, see Figure~1.
Consider a site $\frac{1}{2}$-percolations with free BC as above.
Let \ts $\cU:=\big\{\exists \ts x\in V \. : \. x\lra 12, \. x\lra 34, \. x\lra 56\big\}$ \ts
be the joint connectivity property of the percolation graph.  It was computed by
Simmons~\cite{Sim} (see also~\cite{FZS}), that
\ts $\bP_{\frac12}(\cU)=0.2556897...$  \ts in the limit \ts $\ell\to \infty$.
Consider a uniform random coloring \ts $f: V\to \{a,b,c,d\}$.  Then \eqref{eq:main2} gives:
$$\bP_{\frac12}(\cU)^2 \. = \. 0.0653772... \. \ge \.
\bP\big(\cU_{ab} \cap \cU_{ac} \cap \cU_{ad}\big) \. \ge \. \bP_{\frac12}(\cU)^3 \. = \. 0.0167162...
$$
in the limit \ts $\ell\to \infty$.  Similarly, the inequality~\eqref{eq:main} gives:
$$
\bP\big(\cU_{ab} \cap \cU_{ac} \cap \cU_{bc}\big) \. \le \. \bP_{\frac12}(\cU)^3 \. = \. 0.0167162...
$$
in the limit \ts $\ell\to \infty$.  For \ts $\ell=30$, the sampling of \ts $N=64000$ \ts trials gives
\. $\bP\big(\cU_{ab} \cap \cU_{ac} \cap \cU_{ad}\big)= 0.0172\pm 0.0005$ \. and \.
$\bP\big(\cU_{ab} \cap \cU_{ac} \cap \cU_{bc}\big)= 0.0166\pm 0.0005$.
We conjecture that both probabilities are $\ts \bP_{\frac12}(\cU)^3 = 0.0167162...$ \ts in the limit \ts $\ell\to \infty$.

\smallskip

\subsection*{New critical probability}
Let \ts $G=(V,E)$ \ts be an infinite connected graph.
Consider a uniform random coloring \ts $f: E\to \{a,b,c,d\}$.
For a vertex $x\in V$, consider
\begin{equation}\label{eq:inf-conn}
P(x) \ := \ \bP\big(x\lra_{ab} \infty, \. x\lra_{ac} \infty, \. x\lra_{ad} \infty\big),
\end{equation}
where \ts $x\lra_{st} \infty$ \ts means that $x$ belongs to an infinite cluster
of $st$-colored edges.  Now \eqref{eq:main2} gives:
\begin{equation}\label{eq:inf-bound}
\bP_{\frac12}\big(x\lra \infty)^2 \, \ge \, P(x) \, \ge \, \bP_{\frac12}\big(x\lra \infty)^3.
\end{equation}

Suppose \ts $G=(V,E)$ \ts is a lattice with critical probability
\ts $p_c<\frac12\ts.$ For \ts $\al \in [0,\frac14]$,
consider a random $5$-coloring \ts $f: E\to \{a,b,c,d,\diamond\}$, where
the probabilities of colors $a,b,c,d$ are~$\al$, and the probability of~$\ts\diamond\ts$ is $\ts (1-4\al)$.
Then \ts $E_{ab}$, \ts $E_{ac}$ \ts and \ts $E_{ad}$ \ts are pairwise independent $2\al$-percolations.
Denote by \ts $P_\al=P_\al(x)$ \ts the probability given by \eqref{eq:inf-conn} in this deformation.
Define the following \emph{critical probability} \ts for the colored percolation:
$$\al_c \. := \. \sup\big\{\al \. : \. P_\al(x)=0\big\}.
$$
Now \eqref{eq:inf-bound} implies that \ts $\al_c \le \frac{1}{2} \. p_c$ \ts while
the examples above suggest \ts $\al_c = \frac{1}{2} \. p_c\ts.$  The
numerical experiments also seem to confirm this.  We tested the colored bond 
and site percolations on a triangular lattice with \. $p_c = 2\sin \frac{\pi}{18} = 0.3473...$ \ts and
\. $p_c=\frac12$\ts, respectively \cite{SE1}.  Similarly, we tested the colored bond and site percolations
on a cubic lattice $G=\zz^3$  with \. $p_c=0.2488...$ \ts and \ts $p_c=0.3116...$, respectively
\ts (see e.g.~\cite{SA}).
The results are given in Figure~2.

\smallskip

\begin{figure*}[hbt]
\begin{center}
	\includegraphics[height=4.2cm]{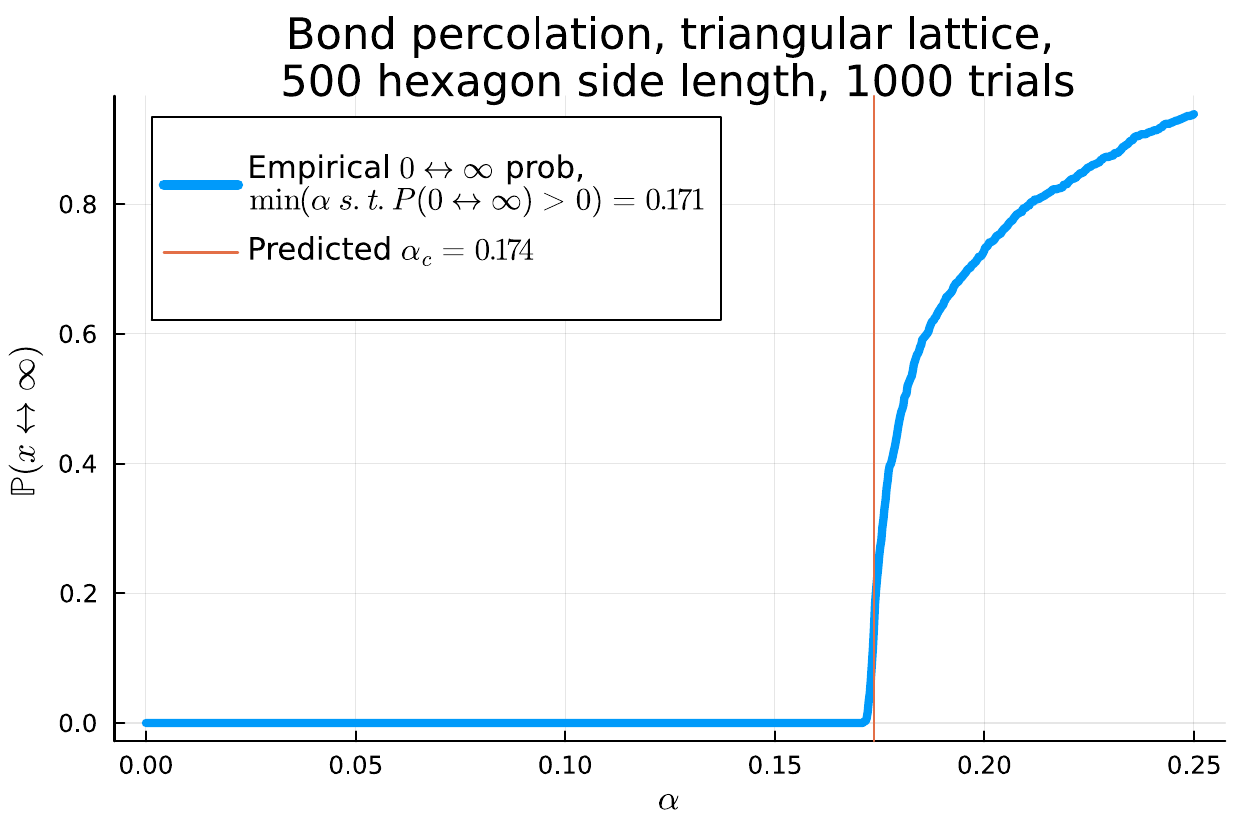}
\hskip1.cm
	\includegraphics[height=4.2cm]{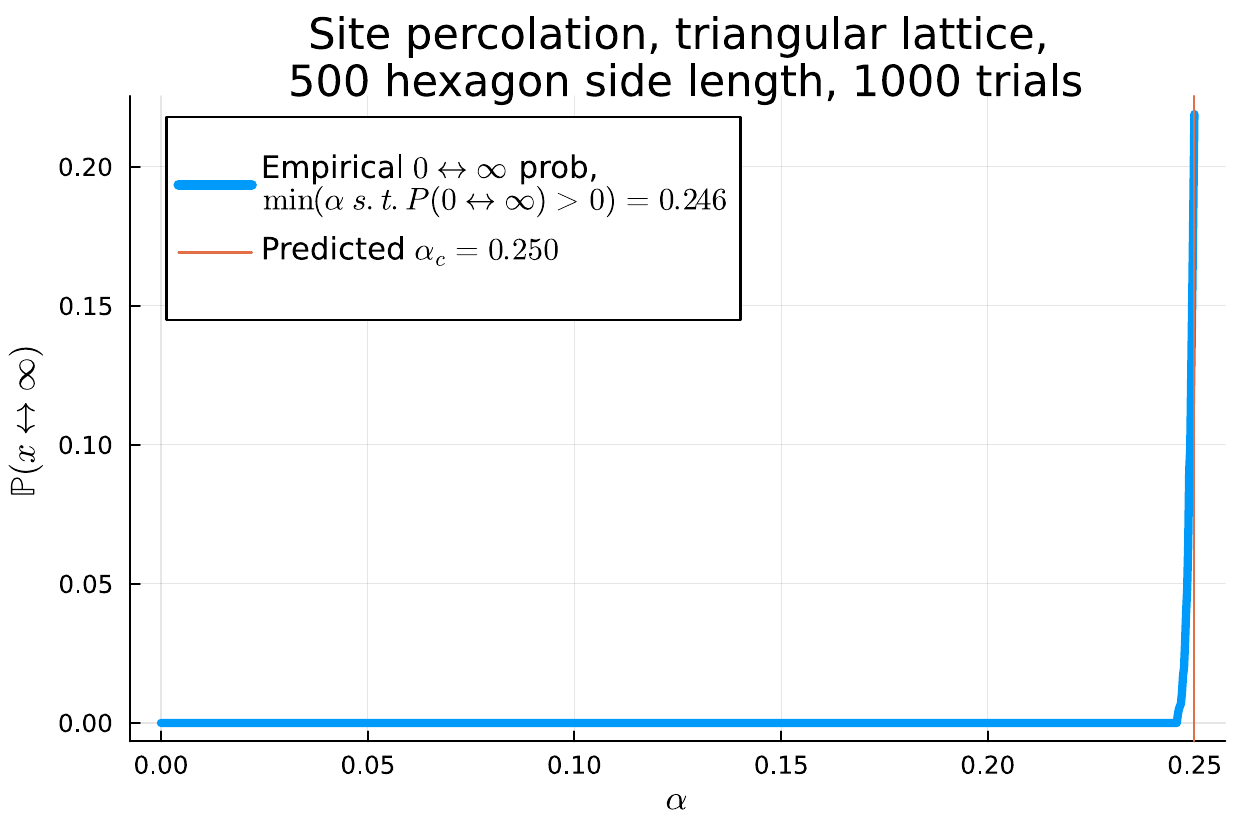}
\end{center}
\begin{center}
	\includegraphics[height=4.2cm]{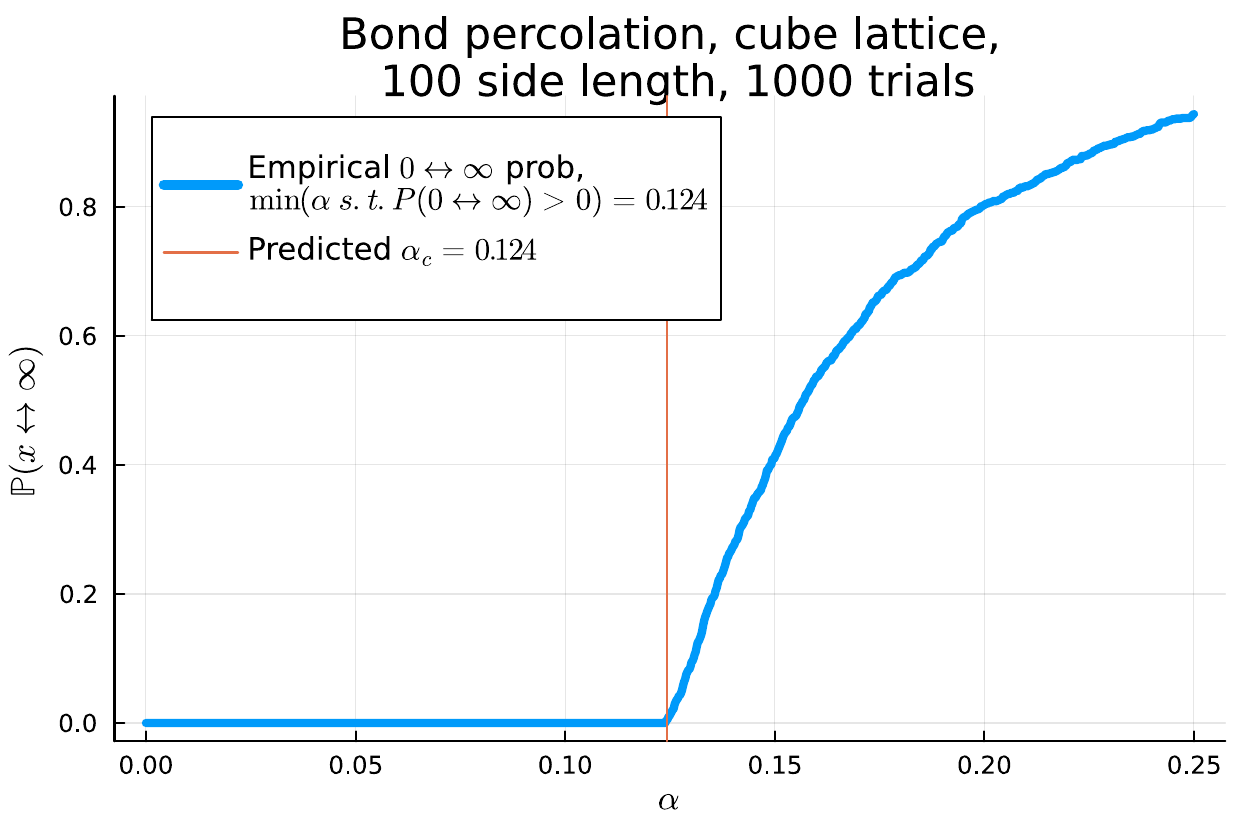}
\hskip1.cm
	\includegraphics[height=4.2cm]{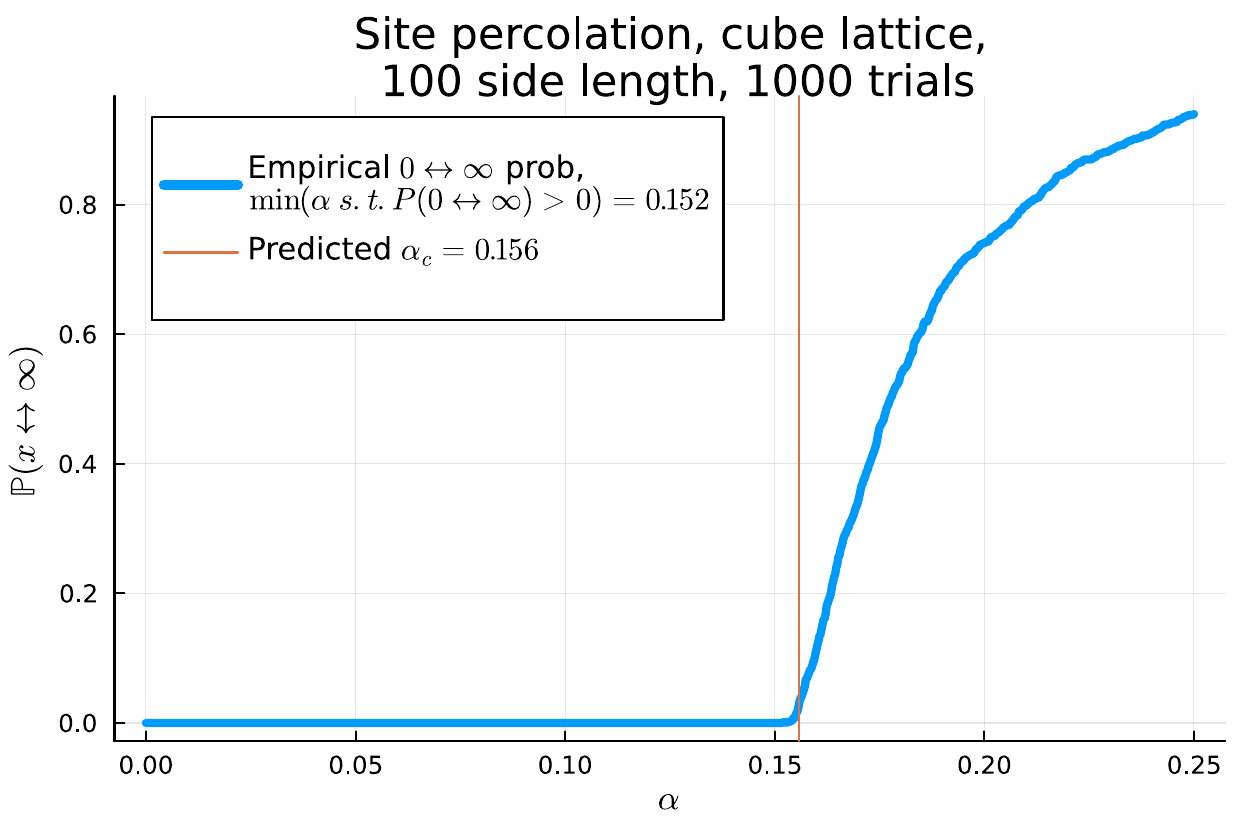}
\caption{Colored bond/site percolations in triangular and cubic lattices.}
\end{center}
\label{f:crit}
\end{figure*}

\subsection*{Conclusions}
The subject of positive dependence for colored percolation is largely unexplored
and can be viewed as a special case of algebraic inequalities for cumulants of
positive functions.  The latter has been actively studied (see \cite{Gla,LS}
for recent references), but the type of inequalities we consider are new.

In full generality, our results extend the Harris--Kleitman inequality~\eqref{eq:HK}
to multiple pairwise independent events.  Following classical approach, this allows
us to give lower and upper bounds on the mutual dependence of these events, and
 to define critical constants $\al_c$ for a deformation of the colored percolation.
Our lower and upper bounds are asymptotically tight for the conjectured crossing
probabilities of the colored percolation on lattices, exhibiting the same phenomenon
as the majority property.

\smallskip

\subsection*{Acknowledgements}
We are grateful to Matija Buci\'c, Tom Hutchcroft, Jeff Kahn and
Bhargav Narayanan for interesting conversations on the subject.
Special thanks to Geoffrey Grimmett and Bob Ziff for helpful remarks
on the paper, and to Aleksandr Zimin for help with numerical
experiments.  The second author was partially supported by the NSF.

\end{document}